\documentclass[times,review]{elsarticle}
\usepackage{pdflscape}
\usepackage{graphicx}
\usepackage{amssymb}
\usepackage{amsmath}
\usepackage{hyperref}
\usepackage{lineno}
\usepackage{subfigure}
\usepackage{tikz}
\usepackage{adjustbox}
\usepackage{soul}
\usepackage{todonotes}
\usetikzlibrary{shapes.geometric, arrows}
\usepackage{catchfilebetweentags}
\usepackage{booktabs}
\usepackage{framed,multirow}

\newcommand{\ffrac}[2]{\ensuremath{\frac{\displaystyle #1}{\displaystyle #2}}}

\UseRawInputEncoding

\modulolinenumbers[5]

\begin{document}


\begin{frontmatter}

    \title{a Decision-Tree based Moment-of-Fluid (DTMOF) Method in 3D rectangular
        hexahedrons}%
    \author[1,2]{Zhouteng Ye}
    \ead{yzt9zju@gmail.com}
    
    \author[2]{Mark Sussman}
    \ead{sussman@math.fsu.edu}
    
    \author[1]{Yi Zhan}
    \ead{yi.zhan@zju.edu.cn}
    
    \author[1]{Xizeng Zhao\corref{cor1}}
    \cortext[cor1]{Corresponding author:
        email: xizengzhao@zju.edu.cn;
    }

    \address[1]{
        Ocean College, Zhejiang University, Zhoushan 316021, Zhejiang, People’s Republic of China}
    \address[2]{
        Department of Applied and Computational Mathematics, Florida State University, United States}


    \begin{abstract}
        The moment-of-fluid (MOF)  method is an extension of the volume-of-fluid method with piecewise linear interface construction (VOF-PLIC).
        By minimizing the least square error of the centroid of the cutting polyhedron,
        the MOF method reconstructs the linear interface without using any neighboring information.
        Traditional MOF involves iteration while finding the optimized linear reconstruction.
        Here, we propose an alternative approach based on a machine learning algorithm: Decision Tree algorithm.
        A training data set is generated from a list of random cuts of a unit cube by plane.
        The Decision Tree algorithm extracts the input-output relationship from the training data,
        so that the resulting function determines the normal vector of the reconstruction plane directly,
        without any iteration.
        The present method is tested on a range of popular interface advection test problems.
        Numerical results show that our approach is much faster than the iteration-based MOF method
        while provides compatible accuracy with the conventional MOF method.

    \end{abstract}


\end{frontmatter}



\section{Introduction}
\label{sec:introduction}

A lot of scientific and engineering problems involve tracking the interface between different materials.
Multiple volumes tracking/capturing methods,
such as volume-of-fluid (VOF) method \cite{hirt_volume_1981,youngs_time-dependent_1982,zhang_new_2008},
level set method \cite{osher_fronts_1988,sussman_level_1994,osher_level_2003},
and front tracking method \cite{unverdi_front-tracking_1992,tryggvason_front-tracking_2001}
are introduced to describe the motion of the interface explicitly or implicitly.
Among those methods, the volume-of-fluid method with piece-wise line interface construction (VOF-PLIC) is one
of the most widely used methods in tracking the interface within the Eulerian framework.

In VOF-PLIC method, the material interface of the material is described with the volume fraction.
When discredited with 3D rectangular hexahedron, the volume fraction can be expressed as
\begin{equation}
    \label{eq:voffrac}
    C_{i,j,k}=\left\{\begin{array}{ll}
        \ffrac{\iiint_{\Omega_{i,j,k}} f(x,y,z) \text{d} x \text{d} y \text{d} z }{\Delta x \Delta y \Delta z} ,
           & \text { Interface cell}     \\ [10pt]
        1, & \text { Material cell}      \\  [6pt]
        0, & \text { Non-material cell }\end{array}\right. .
\end{equation}
Where ${\Omega_{i, j,k}=\left[x_{i-1 / 2}, x_{i+1 / 2}\right] \times\left[y_{j-1 / 2}, y_{j+1 / 2}\right] \times\left[z_{k-1 / 2}, z_{k+1 / 2}\right]}$
is cell domain and $C_{i,j,k}$ is the volume fraction of the color function $f(x,y,z)$ within the cell domain $\Omega_{i,j,k}$. 
Conventional VOF-PLIC method reconstructs the normal vector of the reconstructed interface
by using the a stencil that contains the information of the neighboring grids,
for example, Parker and Youngs' algorithm \cite{parker_two_1992},
mixed Youngs-centered algorithm (MYC) \cite{aulisa_interface_2007},
and the efficient least squares volume-of-fluid interface reconstruction algorithm (ELVIRA) \cite{pilliod_second-order_2004}.
Although some of the VOF-PLIC reconstruction algorithms are second-order accuracy,
when there is not enough information from the neighboring grid,
for example, very small scale droplets,
VOF-PLIC algorithm may not reconstruct the interface accuracy.

Moment of Fluid (MOF) method \cite{dyadechko_moment--fluid_2005,dyadechko_reconstruction_2008}
provides an alternative way to determine the normal vector.
In MOF method, both of centroid $\mathbf{c}= \{c_x,c_y,c_z\}$ and the volume fraction $C$
are used to determine the normal vector of the reconstruction plane $\mathbf{n \cdot x} = \alpha$.
Without using data from adjacent cells,
MOF reconstruction resolves the interface with a smaller minimum scale than the VOF-PLIC algorithm
and has been extended from Cartesian grid to multiple frameworks such as adaptive mesh refinement(AMR)\cite{ahn_multi-material_2007, jemison_coupled_2013,liu_moment--fluid_2020},
arbitrary Lagrangian-Eulerian (ALE) \cite{galera_2d_2011,breil_multi-material_2013}.
It is easy to determine the centroid and volume fraction from the given plane (referred to as forward algorithm $\mathcal{F}$)
\begin{equation}
    \label{eq:forwrd}
    ({\mathbf{c}}, C) = \mathcal{F}({\mathbf{n}}, \alpha).
\end{equation}
Unfortunately, find the cutting plane from the centroid and volume fraction  (referred to as backward algorithm $\mathcal{G}$)
\begin{equation}
    \label{eq:backward}
    ({\mathbf{n}}, \alpha) = \mathcal{G}({\mathbf{c}}, C)
\end{equation}
is not as simple as the forward algorithm.
Eq. \eqref{eq:backward} is typically solved with an iteration algorithm that minimizes
the $L_2$ norm between the reconstructed centroid
and the reference centroid.
The iteration algorithm starts with an initial guess of the normal vector.
At each of the iteration step,
the volume fraction, centroid and the gradient of the objective function are calculated
and used to determine the normal vector for the next iteration step.
In most of the MOF algorithm, the forward algorithm $\mathcal F$ in Eq. \eqref{eq:forwrd} is solved
at every iteration step of backward algorithm $\mathcal G$.

The original MOF algorithm by \citet{dyadechko_moment--fluid_2005} is time-consuming because
a complex polyhedra intersection algorithm is used as the forward algorithm $\mathcal{F}$ to solve Eq. \eqref{eq:forwrd},
and the forward algorithm $\mathcal F$ has to been used 5 times at each iteration to determine the gradient of the
objective function.
Several approaches have been used to accelerate the MOF reconstruction.
\citet{jemison_coupled_2013} proposed a coupled level-set and moment-of-fluid (CLSMOF) by coupling the level set function with MOF.
The level set function is used to provide a better initial guess of the normal vector so that the iteration with fewer steps.
\citet{chen_improved_2016} developed an analytic gradient for the objective function of the MOF iteration.
By using an analytic gradient form,
the number of calling the forward algorithm $\mathcal{F}$ reduced from 5 to 1 time at each iteration.
The algorithm is found to be 3-4 times faster than the original MOF by \citet{dyadechko_moment--fluid_2005}.
Besides boosting the iteration algorithm,
\citet{lemoine_moment--fluid_2017} made their first attempt to derive an analytic form of that describes
Eq. \eqref{eq:backward} as the minimum distance from the reference centroid to a closed, continuous curve.
This is a fully analytic 2D MOF algorithm as a solution to Eq. \eqref{eq:backward} can be obtained by computing the cubic or quartic roots of polynomials instead of iteration.
Unfortunately, this approach cannot be extended to 3D.
\citet{milcent_moment--fluid_2020} proposed an analytic approach to determine the objective function and its gradient instead of the geometrical approach.
Although analytic gradient is much more efficient than the numerical gradient algorithm by \citet{dyadechko_moment--fluid_2005,chen_improved_2016},
iteration is still unavoidable while solving Eq. \eqref{eq:backward}.

The machine learning technique provides a new approach to model the non-linear
input-output function.
It constructs the input-output function by algorithmic learning of
essential features in the training data-set,
rather than deriving the functional relationship using some physical assumption or analytic relationship.
In recent years, machine learning technique has been used in modeling multiphase flow,
and has shown its potential in boosting the performance of the numerical simulation.
For example, \citet{ma_using_2015,ma_using_2016} use neural networks algorithm to enclosure the unknown
terms in average flow.
\citet{qi_computing_2019} estimate the curvature of the VOF-PLIC method by
using the volume fraction of the surrounding cells.
This new approach of estimating curvature has been extended to different frameworks,
such as CLSVOF method \citep{haghshenas_curvature_2019}, level-set method \cite{cardenas_deep_2020}.
\citet{ataei_nplic_2020} proposed a model trained from a data-set of PLIC solutions,
the result shows that the data-driven approach maintains the accuracy of PLIC method at a fraction
of the usual computational cost.
A discussion in the context of multi-phase flow and machine learning algorithm can be found in \citet{gibou_sharp_2018}.

In this study, we apply a machine learning algorithm, called Decision Tree (DT) algorithm
to model the normal vector of the reconstruction plane from
the volume fraction and the centroid in one cell.
The new MOF method is called DTMOF (Decision Tree boosted Moment of Fluid).
The main objective of our DTMOF method is to build an efficient MOF reconstruction function 
for practical multi-phase simulation.
A synthetics data-set is generated from a list of linear reconstruction data.
The resulting functional relationship for MOF reconstruction determines the optimal normal vector directly, 
without any iteration.
The decision tree models the normal vector of the reconstruction plane from the training data.
Our DTMOF model is tested with static reconstruction and several advection cases.
The layout of the paper is as follows:
Section 2 introduces our DTMOF method,
The static reconstruction is tested in Section 3 and compared with other machine learning algorithms.
Several advection cases are tested in Section 4 and finally the conclusion is drawn in Section 5.

It should note that the run-time ratio and robustness of the method could be implementation-dependent.
Out implementation of the code and test cases are available on our Github repository (https://github.com/zhoutengye/NNMOF).
All the cases are done on a workstation with Intel(R) Xeon(R) Platinum 8270 processors with the Intel Fortran compiler 2020 on Linux Mint 19.3.

\section{Decision Tree boosted Moment of Fluid Method}
\label{sec:dtmof}
\subsection{Revisit to Moment-of-fluid reconstruction}
\label{sec:dtmof-mof}
In fluid simulation with MOF method,
the known reference centroid $\mathbf{c}_{\rm{ref}}$ and volume fraction $C_{\rm{ref}}$
may not simultaneously satisfy with a linear cut-off.
To keep the volume conservation,
the MOF algorithm sacrifices the exact centroid matching and looks
for a linear cut-off with the given volume fraction
which provides the best approximation to the reference centroid.

The linear cut-off plane in a 3D rectangular hexahedron cell is defined as
\begin{equation}
    \label{eq:mofplane}
    \mathcal{B} = \left\{\mathbf{x} \in \mathbb{R}^{3} \mid \mathbf{n} \cdot\left(\mathbf{x}-\mathbf{x}_{0}\right)+\alpha=0\right\},
\end{equation}
where $\mathbf{n}$ is the normal vector, $\mathbf{x_{0}}$ is the reference point of the cell,
either the center of the cell or the lower corner of the cell,
depending on the computational algorithm.
$\alpha$ is the parameter that represents the distances from the reference point $\mathbf{x_0}$.
The volume fraction of the reconstruction polyhedron $C_A$ should be equal to the reference volume fraction
\begin{equation}
    \label{eq:mofvolumeconstraint}
    \left|{C}_{\mathrm{ref}}(\mathbf{n}, \alpha)
    -{C}_{A}(\mathbf{n}, \alpha)\right|=0.
\end{equation}
In addition to the constraint on volume fraction,
the MOF reconstruction also minimizes error of the centroid
\begin{equation}
    \label{eq:mofcentroidconstraint}
    E_{\mathrm{MOF}}=\left\|\mathbf{c}_{\mathrm{ref}}-\mathbf{c}_{A}(\mathbf{n}, \alpha)\right\|_{2}.
\end{equation}
The normal vector can either be represented with the vector form $\mathbf{n}=(n_x,n_y,n_z)$
or spherical coordinate form $\mathbf{\Phi} = (\phi,\theta)$.
The conversion between the two forms are 
\begin{equation}
    \label{eq:norm2angle}
    \mathbf{n(\phi, \theta)}=
    \left(\begin{array}{c}
            \sin (\phi) \cos (\theta) \\
            \sin (\phi) \sin (\theta) \\
            \cos (\phi)
        \end{array}\right),
\end{equation}
\begin{equation}
    \label{eq:angle2norm}
    \mathbf{\Phi}(n_x, n_y, n_z)=
    \left(\begin{array}{c}
            \arctan(\frac{n_y}{n_x}) \\
            \arctan(\frac{\sqrt{n_x^2+n_y^2}}{n_z}).
        \end{array}\right).
\end{equation}
With the constraint of the volume fraction in Eq. \eqref{eq:mofvolumeconstraint},
$\alpha$ can be uniquely defined by the known normal vector.
Substitute Eq. \eqref{eq:norm2angle} into Eq. \eqref{eq:mofcentroidconstraint}
and the objective function of the centroid is simplified as a function of $\phi$
and $\theta$. Minimizing the error $E_{\rm{MOF}}$ is to find
the optimized $(\phi^{*}, \theta^{*})$
\begin{equation}
    \label{eq:mofminimizeangle}
    E_{\mathrm{MOF}}\left(\phi^{*}, \theta^{*}\right)=\left\|\mathbf{f}\left(\phi^{*}, \theta^{*}\right)\right\|_{2}=\min _{(\phi, \theta): {\mathrm{Eq}}. (4) \text { holds }}\|\mathbf{f}(\phi, \theta)\|_{2}
\end{equation}
where,
\begin{equation}
    \label{eq:parafunc}
    \mathbf{f}: \mathbb{R}^{2} \rightarrow \mathbb{R}^{3}, \quad \mathbf{f}(\phi, \theta)=\left(\mathbf{c}_{\text {ref }}-\mathbf{c}_{A}(\phi, \theta)\right)
\end{equation}
The minimization problem in Eq. \eqref{eq:mofminimizeangle} is a non-linear least square problem for $\phi$ and $\theta$,
which is solved numerically with an optimization algorithm.

\ExecuteMetaData[figures.tex]{fig:sketchmof}

A typical solution procedure for an iteration-based optimization algorithm is

0. Choose initial angles $(\theta_0$ , $\phi_0)$ and initialize iteration step $k = 0$.

While not converged,

1. Find $\alpha_k (\phi_k , \theta_ k)$ such that Eq. \eqref{eq:mofvolumeconstraint} holds.

2. Find the centroid $\mathbf{c_k}$ $(\alpha_k$ , $\phi_k$ , $\theta_k)$.

3. Estimate the gradient of the objective function $(\frac{\partial \mathbf{f}}{\partial \phi}_k,\frac{\partial \mathbf{f}}{\partial \theta}_k$).

4. Update the angles: $\left(\phi_{k+1}, \theta_{k+1}\right)$

5. $k:=k+1$

The above procedure can be estimated using non-linear optimization method,
for example,
BFGS algorithm \citep{ahn_multi-material_2007, chen_improved_2016,milcent_moment--fluid_2020}
or Gauss-Newton \citep{jemison_coupled_2013,asuri_mukundan_3d_2020} algorithm.
Estimating the gradient of the objective function
$(\frac{\partial \mathbf{f}}{\partial \phi}_k,\frac{\partial \mathbf{f}}{\partial \theta}_k$)
takes most of the computational time during the iteration.
Conventional MOF method \citep{ahn_multi-material_2007} estimates the gradient numerically with a
computational geometrical algorithm.
\citet{chen_improved_2016} proposed an analytic geometric algorithm for
$(\frac{\partial \mathbf{f}}{\partial \phi}_k,\frac{\partial \mathbf{f}}{\partial \theta}_k$)
which reduce the gradient estimation from 5 to 1 time in each iteration step.
\citet{milcent_moment--fluid_2020} proposed an analytic form for the
numerical gradient
$(\frac{\partial \mathbf{f}}{\partial \phi}_k,\frac{\partial \mathbf{f}}{\partial \theta}_k$) ,
which significantly reduce the computational cost of the gradient estimation.
Although the analytic MOF algorithm of \citet{milcent_moment--fluid_2020}
is much more efficient than the conventional MOF algorithm,
iteration is unavoidable in order to get the optimized angle.

\subsection{Machine learning approach for MOF reconstruction}
In this study,
a machine learning approach is used to determine the
optimized cut-off of the MOF algorithm.
The key steps for the machine learning algorithm are:
\begin{itemize}
    \item Generating two synthetic data-set from the linear cut-off,
          one for training and the other for test.
    \item Fitting the training data set using machine learning algorithm (the learning stage) to
          find the function that determines the optimized angle ${\mathbf{\Phi}}$
          from the reference centroid $\mathbf{c}_{\rm{ref}}$
          and the reference volume fractions $C_{\rm{ref}}$.
    \item Testing the function on test data set.
    \item Predict the angle $\tilde{\mathbf{\Phi}}$ using the functional relationship built from training.
\end{itemize}
Note that in this section, we use the symbol \textasciitilde \hspace{1mm}
to distinguish the prediction value from the data value.

The training and test data sets are generated from
a list of random linear cut-offs.
A machine learning algorithm: Decision Tree algorithm
(See next subsection) is used to build the functional
relationship between the known centroid-volume fraction pair
$(\mathbf{c}_{\rm{ref}},C_{\rm{ref}})$
and the optimized angle $\tilde{\mathbf{\Phi}}$.
We do not build the functional relationship between
$(\mathbf{c_{\rm{ref}}},C_{\rm{ref}})$ and
$\tilde {\mathbf{\Phi}}$ directly.
Instead, we build the relationship as a guess-correction procedure.
The initial guess of the normal vector
\citep{ahn_multi-material_2007,milcent_moment--fluid_2020}
corresponds with the normal vector from
the centroid towards the center and the
normal vector is parameterized with Eq. \eqref{eq:norm2angle}
\begin{equation}
    \label{eq:initialguess}
    \mathbf{\Phi_0}=  \mathbf{\Phi(c_{\rm{ref}}-x^0)},
\end{equation}
where $\mathbf{x^0}$ is the center of the hexahedron grid.
$\mathbf{\Phi}_0$ is used as the input feature instead of the centroid $\mathbf{c_{ref}}$.
The correction from the initial guess to the exact angle is expressed as
\begin{equation}
    \label{eq:dtmof2}
    \Delta \mathbf{\Phi} = \mathbf{\Phi} - \mathbf{\Phi_0}.
\end{equation}
The functional relationship in our DTMOF model is now defined as
\begin{equation}
    \label{eq:dtmof3}
    \Delta \mathbf{\Phi}  = f (\mathbf{\Phi_0}, {C})
\end{equation}
The function contains three input features and two output target values.
The goal of the data training is to get a prediction of $\mathbf{\Delta \Phi}$  corresponding to the inputs, $\mathbf {x_c}$ and $C$.
In the prediction session, the value of $\mathbf{\Phi_0}$ is firstly determined from
Eq. \eqref{eq:initialguess} and the prediction of the angle $\mathbf{\Phi}$ is made by
\begin{equation}
    \label{eq:deangle}
    \mathbf{\tilde \Phi} = \mathbf{\Phi_0} + f(\mathbf{\Phi_0}, C).
\end{equation}

\ExecuteMetaData[figures.tex]{fig:xyzptf}

It is important for the training data to cover the range of all
possible inputs in fluid simulation.
\citet{dyadechko_moment--fluid_2005}
show that for linear reconstruction plane,
the volume of the cutting polyhedron $C$
is uniquely identified by its centroid $\mathbf{c}$.
For 3D cut-offs from a hexahedron, the locus of the centroid of with fixed volume is a closed surface \citep{milcent_moment--fluid_2020}.
We plot the loci of the centroids for fixed volume fraction in $\Omega_2(c_x,c_y,c_z)$
in Fig. \ref{fig:xyzptf} (a).
Those closed surfaces can be mapped to the data space
$\Omega(\phi, \theta, C)$ as planes (See Fig. \ref{fig:xyzptf} (b)).
When the training data is generated from uniform distribution 
data space $\Omega_2$, 
it ensures the coverage of all possible loci of centroids 
in the unit cube $\Omega_1$.
For an arbitrary hexahedron $\Omega_h$ with variable edge length,
a mapping from $\Omega_h$ to the unit cube $\Omega_1$ could be done 
to find the optimized angle.
Similar mapping can also be found in the data-driven MOF approach of \citet{cutforth2021efficient}.

\subsection{Decision Tree algorithm}
\label{sec:dtmof-dt}

Decision Tree (DT) method is a machine learning used for
classification and regression problem \citep{madsen_methods_2004}.
The idea of the machine learning method is to find patterns between input features and target values through the data training process.

\ExecuteMetaData[figures.tex]{fig:dt}

The training data $D$,
including input features and output targets,
are generated using a list of random polyhedra of a plane cutting the unit cube.
Detailed data generation is described in the next section.
During the training phase,
the DT algorithm splits the data-set into two smaller partitions (branches) recursively,
as shown in Fig. \ref{fig:dt}.
The final tree structure is determined by adjusting the algorithm recursively
with the objective to minimize the sum of variances in the
response values across all the partitions (leaves).
The best split of the subset $X \subset D$,
decision tree algorithm splits the data into two smaller subsets
\begin{equation}
    \label{eq:split}
    R_{1}(j, s)=\left\{X | X_{j} \leq s\right\} \text { and } R_{2}(j, s)=\left\{X | X_{j}>s\right\},
\end{equation}
where $j$ means $j$-th input feature,
$s$ indicates the value of the threshold
$R_{1}(j,s) y\cup R_{2}(j,s) = X$.
The best split $(s_{\rm{best}})$ of the subset $X$ can be determined by
finding the minimal mean square root (MSE) of the output value variable $y$ across all
possible splits.
\begin{equation}
    \label{eq:mse}
    E_{s_{\rm{best}}} = \rm{min} \left [\sum_{i: x_{i} \in R_{1}(j, s)}\left(y_{i}-\hat{y}_{R_{1}}\right)^{2}+\sum_{i: x_{i} \in R_{2}(j, s)}\left(y_{i}-\hat{y}_{R_{2}}\right)^{2} \right ].
\end{equation}
After all subsets reaching certain criteria (e.g. max depth level, maximum impunity),
the training finishes and the functional relationship between
in Eq. \eqref{eq:dtmof3} is stored in the decision tree.

Note that \citet{cutforth2021efficient} also applied a data-driven approach to accelerate the MOF reconstruction in 2D.
In \citet{cutforth2021efficient} all training data are used as the database during computing.
For a $800^2$ training data-set from 2D grid, it takes about 1.3 GB disk space,
when extending to 3D, the size of the data would be challenging for the computer, 
especially for parallel computing on distributed memory.
The DT algorithm splits the training data into smaller subsets of data on its
leaf nodes, averaging the value of the data in each subset.
So that the DT prediction uses much less disk space than the training data.

\section{Data generation and training}
\label{sec:trainig}

In this section,
we generate two synthetic data-sets,
one for training and the other for test.
We also compare accuracy and efficiency of our DTMOF algorithm
with the original MOF algorithm and other two machine learning algorithms:
Neural Networks algorithm and Random Forest algorithm.

The training and test data sets are both generated from a list of
planes cutting a unit cube.
The initial guess of the angle $\mathbf{\Phi_0}$
and the volume fraction of the cut-off polyhedron $C$ are used as the input features
and the correction of the angle $\mathbf{\Delta \Phi}$ is used
as the output target.
The training data-set contains $1\times 10^9$ sets of data and
test data sets contains $1\times10^8$ sets of data.
In this study, we compare the efficiency and accuracy of the decision tree
with other machine learning algorithms.
We use two criteria to estimate the training and prediction:
Mean $L_1$ error of centroid $E_c$ and coefficient of determination $R^2$

\begin{equation}
    \label{eq:error_centroid}
    E_{c}=
    \frac{\sum_{i, j}\left| \mathbf{\tilde{c}} - \mathbf{c}_A\right|}{N},
\end{equation}

\begin{equation}
    \label{eq:centroid_error}
    R^2= 1 - \frac{1}{2}\left (
    \frac{\sum_{i}^{N}\left(\theta_{i}^e-\Tilde{\theta}_{i}\right)^{2}}{\sum_{i}^{N}\left(\theta_{i}^e-\bar{\theta}\right)^{2}} +
    \frac{\sum_{i}^{N}\left(\phi_{i}^e-\Tilde{\phi}_{i}\right)^{2}}{\sum_{i}^{N}\left(\phi_{i}^e-\bar{\phi}\right)^{2}}
    \right ),
\end{equation}
where $N$ is the size of data set, $\theta^e$ and $\phi^e$ are exact value of the angles.

\ExecuteMetaData[figures.tex]{fig:DTConvergence}

\ExecuteMetaData[figures.tex]{fig:dperror}

\ExecuteMetaData[figures.tex]{fig:dterror}

The maximum depth of the decision tree $d_{max}$ plays an important role in DTMOF correction.
When there is no limit of the maximum depth of the tree,
the amount of the tree nodes would be huge,
which would affect both computational efficiency and storage.
While if the maximum depth is too small,
the decision tree may not be able to model the input-output relationship correctly.
The CPU time, total leaf count and $R^2$ with the change of maximum layer are plotted in Fig. \ref{fig:DTconvergence}.
With a deeper decision tree,
the computational cost and storage increase,
which brings in better accuracy.
When $d_{max}<=20$,
the computational cost and storage increases slowly,
while the value of $R^2$ increases rapidly.
When $d_{max}>20$,
$R^2$ changes slowly with the increase of the maximum depth $d_{max}$.
For a balance between among the storage, CPU time and accuracy,
we choose $d_{max}=20$ in this study.
The decision tree occupies about 11 MB disk space,
the $R^2$ value for training and test data sets are 0.977 and 0.990, respectively.

The isosurfaces of the two output variables ($\Delta \phi, \Delta \theta$)
in the data space $\Omega_2$
are plotted in Fig. \ref{fig:dperror} and Fig. \ref{fig:dterror}.
The predicted $\Delta \phi$ and $\Delta \theta$ are compared with the
exact values in the test data set.
Although the isosurfaces is not as smooth as the isosurfaces in the test data set,
the predicted values from the DTMOF algorithm shows an overall
good agreement with the exact values.

\ExecuteMetaData[tables.tex]{tab:mlstatic}

We also compare the DT results with
conventional MOF algorithm,
analytic MOF algorithm,
and two other machine learning algorithms:
Neural Network algorithm and Random Forest algorithm
The results are shown in Table \ref{tab:mlstatic}.
In conventional MOF \citep{dyadechko_moment--fluid_2005} and analytic MOF algorithm \citep{milcent_moment--fluid_2020},
the tolerance of the objective function is $10^{-8}$
and the maximum iteration step is 100.

The multi-layer Neural Network algorithm \citep{goodfellow_deep_2016}
is one of the most popular
machine learning algorithms.
After $1 \times 10^8$ sets of the combinations of hyper-parameters
using GridSearchCV \citep{pedregosa_scikit-learn_2011},
we realize that the accuracy of the artificial neural network mostly relies on
the number of neurons in this problem.
We only select two configurations of the artificial neural network in Table \ref{tab:mlstatic}.
With a small amount of neurons,
the neural network cannot represent the input-output relationship correctly.
With the increasing of neurons,
although the model predicts the correction angles more accurately,
the computational cost increases significantly.

The Random Forest algorithm \citep{tin_kam_ho_random_1995} is an ensembled algorithm which uses
multiple regression trees as estimators.
Although the random forest algorithm is reported to have better performance
especially on preventing the overfitting on a single regression algorithm.
The ensembled regressor takes much more computational cost than a single
regressor.
In this study,
as the distribution of the distribution of data in $\Omega_2$ is a uniformly distribution,
the decision tree algorithm has got as good result as the random forest algorithm.

It should note that, when compared with the original iteration algorithm,
the error of the centroids from the DT prediction is  4-order larger.
However,
in fluid simulation,
the exact reconstruction may not be the linear cut-off,
the optimized linear cut-off results in
the small difference between the optimized centroid and reference.
We show in the next section that in the practical problem,
our DT based algorithm reconstructs the normal vector with satisfying results.

\section{Numerical tests}
\label{sec:numerical}
In this section,
we test the accuracy and efficiency of our proposed MOF method
with some test cases.
The reconstruction algorithm is applied to a 3D advection equation
\begin{equation}
    \begin{aligned}
        \label{eq:advection_eq}
         & \frac{\partial C}{\partial t} + \mathbf{u} \cdot \nabla C = 0, \\
         & \frac{\partial \mathbf{c}}{\partial t} = \mathbf{u}.
    \end{aligned}
\end{equation}
A directional-splitting algorithm applied to solve Eq. \eqref{eq:advection_eq}.
The implementation of the advection algorithm follows \citet{jemison_compressible_2014}.

Three tests are taken in this section,
which represent various different scenarios:
translation, rotation, shear, breaking up, and merging.
We also compare our method with the conventional
MOF method \citep{dyadechko_moment--fluid_2005},
analytic MOF method \citep{milcent_moment--fluid_2020}
and ELVIRA method \citep{pilliod_second-order_2004}.
Again, the implementation of traditional MOF and analytic MOF are
adopted from notes code,
the maximum iteration is 10,
and the tolerance for iteration is $10^{-8}$.

Two different error measurement criteria \citep{zhang_new_2008} are used in this study:

(1) Relative distortion error
\begin{equation}
    \label{eq:relative_error}
    E_{r}=\frac{\sum_{i, j,k}\left|f_{i, j,k}-f_{i, j,k}^{0}\right|}{\sum_{i, j} f_{i, j, k}^{0}},
\end{equation}

(2) Geometrical error
\begin{equation}
    \label{eq:geometrical_error}
    E_{g}=\frac{\sum_{i, j,k}\left|f_{i, j,k}-f_{i, j,k}^{0}\right|}{h^3},
\end{equation}

The order of the accuracy \cite{zhang_new_2008,aulisa_interface_2007} is defined as
\begin{equation}
    \label{eq:accuracy}
    \mathcal{O}_{h}=\log _{2}\left(\frac{E_{g}\left(\frac{1}{2 h}\right)}{E_{g}\left(\frac{1}{h}\right)}\right).
\end{equation}

\subsection{Translation test}
\label{sec:numerical-translation}

Translation test is one of the most basic benchmark tests for interface tracking methods.
We modify the 2D shapes from \citet{rudman_volume-tracking_1997} to 3D
and add an additional shape "letter A" to the translation test.
The initial setup and parameters are shown in Fig. \ref{fig:translationinit}.
With periodic boundary conditions being set up on domain boundaries,
the 4 initial shapes remain unchanged theoretically after one period of evolution
in a uniform constant velocity field.

\ExecuteMetaData[figures.tex]{fig:translationinit}
\ExecuteMetaData[tables.tex]{tab:translationerror}

The ELVIRA method calculates the normal vector using a stencil
that contains the neighboring grids,
which leads to the smear-out of the sharp corners
(See 2D results in Fig. \ref{fig:translationresult2d} and
3D results in Fig. \ref{fig:translationresult3d}).
\ExecuteMetaData[figures.tex]{fig:translationresult2d}

\ExecuteMetaData[figures.tex]{fig:translationresult3d}
The MOF method,
including the conventional MOF,
Analytic MOF and DTMOF,
preserve the sharp corner better compared with ELVIRA result.
In  Fig. \ref{fig:translationresult2d} and Fig. \ref{fig:translationresult3d},
the numerical results of DTMOF and iteration-based MOF do not have
visual difference from each other.

The errors of the interface tracking and run-time ratio are  shown in Table \ref{tab:translationerror}.
The relative error $E_r$ of the DTMOF results are smaller than
the ELVIRA results,
and very close to the analytic MOF
and conventional MOF method.
The run-time ratio shows that the
DTMOF method is about 7 times faster
than the analytic MOF method \citep{milcent_moment--fluid_2020}
and more than 700 times faster than the
conventional MOF method \citep{ahn_multi-material_2007}.
Compared with the static reconstruction test
in Section \ref{sec:trainig} 
the acceleration ratio is smaller.
This is due to the run-time on the
advection algorithm of volume fraction and centroid.

\subsection{Rotation test: Zalesak's disk}
\label{sec:numerical:zalesak}

The Zalesak's disk rotation test is firstly introduced by
\citet{zalesak_fully_1979}
and used by many other studies
\citep{rudman_volume-tracking_1997,aulisa_geometrical_2003,aulisa_interface_2007,zhang_new_2008}.
In Zalesak's disk rotation,
the rotation velocity field is defined with the following stream function
\begin{equation}
    \label{eq:zalesakvelosity}
    \psi(x,y) = -\frac{\omega}{2}[(x-x_0)^2+(y-y_0)^2],
\end{equation}
where $x_0,y_0$ are the center  of the rotation.
\citet{enright_hybrid_2002} modified and extended the problem to 3D
in which the shape is defined with a notched sphere rather than a
notched cube and only rotates.
The third component of the velocity field is set to 0,
and the other two components of the velocity field
remain the same as the 2D problem as defined in Eq. \ref{eq:zalesakvelosity}.
We extend the 2D problem of Zalesak's disk \citet{zalesak_fully_1979}
to 3D in the same way as \citet{enright_hybrid_2002},
the setup of the problem is shown in Fig. \ref{fig:zalesakinit}.
the rotational velocity field $(u_x,u_y)$ is defined
by the stream function Eq. \eqref{eq:zalesakvelosity} with
the value of $\omega$ is $4\pi$.
The velocity component at $z$ direction is a uniform velocity $u_z=0.5$ and periodic boundary
condition is applied at $z$ direction.
After a full revolution of $2\pi$ rotation,
the notched sphere returns to its initial location.

\ExecuteMetaData[figures.tex]{fig:zalesakinit}

\ExecuteMetaData[tables.tex]{tab:zalesakerror}

\ExecuteMetaData[figures.tex]{fig:zalesakresult2d}

\ExecuteMetaData[figures.tex]{fig:zalesakresult3d}

The material interface of the DTMOF results
are compared with the two iteration-based MOF methods and the ELVIRA method
in Fig. \ref{fig:zalesakresult2d} and Fig. \ref{fig:zalesakresult3d}.
The DTMOF methods remains its robustness during the evolution
in the rotation velocity field,
the material interface of the DTMOF method has no visual difference from
the iteration-based MOF methods and better than the ELVIRA method,
especially the sharp corner.

Table \ref{tab:zalesakerror} shows the interface errors are measured by Eq. \eqref{eq:geometrical_error}, the order of the model measured by Eq. \eqref{eq:accuracy} and run-time ratio with respect to the run-time of DTMOF method.
The run-time ratio is the averaged run-time ratio with the three grid resolutions.
The global error shows that the DTMOF method is as accurate as the two iteration-based MOF methods and better than the ELVIRA method.
Although the geometrical of the DTMOF method is slightly larger than the traditional MOF \citep{dyadechko_moment--fluid_2005} and analytic MOF \citep{frank_machine-learning_2020},
however, the difference is negligible.
In this case, the DTMOF method is about 6 times faster than the analytic MOF method \citep{frank_machine-learning_2020} and more than 450 times faster than the conventional MOF method \citep{dyadechko_moment--fluid_2005}.

\subsection{Deformation test: reverse vortex}
\label{sec:numerical-vortex}
The deformation test is firstly introduced in \citet{leveque_high-resolution_1996}
and also used in testing the volume tracking/capturing methods
\citep{enright_hybrid_2002,jemison_coupled_2013,kawano_simple_2016,asuri_mukundan_3d_2020}.
The deformation velocity field is defined as
\begin{equation}
    \label{eq:defoemationvelocity}
    \begin{aligned}
        u_x(x, y, z) & =2 \sin ^{2}(\pi x) \sin (2 \pi y) \sin (2 \pi z) \cos (\pi t/T) \\
        u_y(x, y, z) & =-\sin (2 \pi x) \sin ^{2}(\pi y) \sin (2 \pi z) \cos (\pi t/T)  \\
        u_z(x, y, z) & =-\sin (2 \pi x) \sin (2 \pi y) \sin ^{2}(\pi z) \cos (\pi t/T)
    \end{aligned},
\end{equation}
the flow velocity is time-dependent and in the interval $0 <t  < T$.

The setup of the model is shown in Fig. \ref{fig:rvinit}. In this study, a sphere with the radius $R=0.15$ is located at $\mathbf{x_{0}}=(0.35,0.35,0.35)$,
and the period $T=3$.
Fig. \ref{fig:reversevortexevolution} shows the evolution of the problem.
With the deformation velocity field,
the initial shape of sphere reaches its maximum deformation
at $t=T/2$,
the flow reverses after afterwards and recoveries
to the initial sphere at time $t=T$.

\ExecuteMetaData[figures.tex]{fig:reversevortexinit}

\ExecuteMetaData[figures.tex]{fig:reversevortexevolution}

\ExecuteMetaData[tables.tex]{tab:reversevortexerror}

\ExecuteMetaData[figures.tex]{fig:reversevortex100result}

\ExecuteMetaData[figures.tex]{fig:reversevortexresult2d}

The material interface of the DTMOF results with grid number of
$100 \times 100 \times 100$ at time $t=T/2$ and $t=T$
are compared with the two iteration-based MOF methods and the ELVIRA methods
in Fig. \ref{fig:reversevortex100result}.
Subject to the two rotating vortices,
the initial sphere starts to stretch and part of the interface thin out to
one grid cell at the time $t=T/2$.
All methods failed to resolve the thin topology exactly,
The interface of the DTMOF result has no visual difference from the
traditional iteration-based MOF results,
and shows less deformed than the ELVIRA results.
When the thin topology at $t=T/2$ is under-resolved,
the shape could not recover to the identical sphere under the revered vortex
as shown in the second row in Fig. \ref{fig:reversevortex100result}.
However,
the DTMOF and iteration-based MOF results recover to the initial spherical shape better
compared with the ELVIRA result.

Fig. \ref{fig:reversevortexresult2d} shows the 2D slice of the sphere at $t=T$ for
grid number $50\times50\times50$ and $100\times100\times100$.
Unlike the translation and Zalesak's disk rotation test,
a visual difference between the iteration-based MOF and DTMOF method is observed
in this test.
This could be caused by the
topological change during the simulation due to the insufficient grid resolution.
The topological change makes the test more severe than other tests.
Nevertheless,
the overall shape are very close to each other,
and better than ELVIRA results.

Table \ref{tab:reversevortexerror} shows the interface errors measured by Eq.
\eqref{eq:geometrical_error}
order of the model measured by \eqref{eq:accuracy} and run-time ratio with respect to the run-time of DTMOF results.
Again, the run-time ratio is the averaged run-time ratio with the three grid resolutions.
Both iteration-based MOF and DTMOF results show smaller geometrical error than the ELVIRA
results.
However, the convergence ratio $O_h(50)$ of the ELVIRA results are greater than the other results from MOF,
while $O_h(100)$ of ELVIRA is smaller.
The geometrical error of DTMOF result is slightly larger than the iteration-based MOF results,
but provides a compatible accuracy with those from the iteration-based MOF results.
In this case, the DTMOF is about 8.6 times faster than the analytic MOF \citep{frank_machine-learning_2020} and more than 550 times faster than the conventional MOF \citep{dyadechko_moment--fluid_2005}.

\section{Conclusions}
\label{sec:conclusion}

The machine learning method provides an alternative way to extract a functional relationship
between the input variables and output targets
when there is no basic expression available
or it is too complicated to get the basic expression.
With a proper choice of the training data sets, training method and training parameters,
the machine learning method can build a reasonable well funcational relationship.

In this study,
the machine learning method is used to find the optimized linear cut-off for MOF method.
A guess-correction procedure is used to represent the functional relationship between
the known centroid, volume fraction and the optimized angle.
The training and test data sets are generated from a list of random cut-off from a unit cube
and the functional relationship for the angle correction is done by a machine learning algorithm: Decision Tree algorithm.

Static reconstruction tests show that our DTMOF method fits the training data with
a satisfactory accuracy.
Compared with other machine learning algorithms (Neural networks and Random Forest algorithms) the iteration-based MOF methods (conventional MOF and analytic MOF methods),
the DTMOF has a balance between accuracy and efficiency.
In the reconstruction test,
our DTMOF method is about 18 times faster than the analytic MOF method and about 3000 times faster than the conventional MOF method.
In several advection tests,
our DTMOF method shows a compatible accuracy to the iteration based MOF methods,
however, is more than 6 times faster than the analytic MOF method
and more than 450 times faster than the conventional MOF method.
The results show that our DTMOF method provides accurate and robust results with a lower computational cost compared with the iteration-based MOF method.

In this study,
all computational grids are cube grids.
we have not tested the reconstruction of our DTMOF method on arbitrary
rectangular with different edge lengths.
It is likely that there are significant opportunities to do so.
We only implement the DTMOF algorithm on rectangular grid,
however, it is possible to extend the machine learning boosted approach to other
grid systems.
Especially for unstructured grid,
in which a more complex grid and cut-offs grometry are involved
and no simple and efficient algorithm (like the analytic MOF on hexahedron grid)
available,
the functional relationship from the machine learning approach could potentially
get a higher accelerating ratio.


\section*{Acknowledgments}
The support provided by  National Science Foundation of China (Grant Nos. 51979245, 51679212),
and China Scholarship Council (CSC) and the during a visit of Zhouteng Ye to Florida State University is acknowledged.

\bibliography{refs}

\end{document}